\documentclass{amsproc}
\usepackage{verbatim,graphics,amsmath,amsfonts,amssymb}
\newtheorem{theorem}{Theorem}[section]
\newtheorem{lemma}[theorem]{Lemma}
\newtheorem{proposition}[theorem]{Proposition}

\newtheorem{definition}[theorem]{Definition}

\newtheorem{remark}[theorem]{Remark}
\def\pagenumber{1}

\begin{document}
\setcounter{page}{\pagenumber}
\newcommand{\T}{\mathbb{T}}
\newcommand{\R}{\mathbb{R}}
\newcommand{\Q}{\mathbb{Q}}
\newcommand{\N}{\mathbb{N}}
\newcommand{\Z}{\mathbb{Z}}
\newcommand{\C}{\mathbb{C}}
\newcommand{\tx}[1]{\quad\mbox{#1}\quad}

\title[Pasting and Reversing over rings]{\mbox{}\\[1cm]On Pasting and Reversing operations over some rings}
\author[P. Acosta-Humanez, A.L. Chuquen \& A.M. Rodr\'{\i}guez]{Primitivo B. Acosta-Humanez$^{a,b}$, Adriana Lorena Chuquen$^c$, \& \'Angela Mariette Rodr\'{\i}guez$^d$}

\maketitle
\vspace{-.5cm}

{\footnotesize
\begin{center}
$a$ Escuela de Caminos, Canales y Puertos - Universidad Polit\'ecnica de Madrid, Spain\\
$b$ Instituto de Matem\'aticas y sus Aplicaciones (IMA) - Universidad Sergio Arboleda, Santa Marta, Colombia\\
$c$ Escuela de Matem\'aticas - Universidad Sergio Arboleda, Bogot\'a, Colombia\\
$d$ Facultad de Ciencias, Departamento de F\'{\i}sica - Universidad Nacional de Colombia, Bogot\'a
\end{center}}
\vspace{.5cm}
\centerline{\footnotesize \it To Rodrigo Noguera Barreneche, 1892--1972.}
\bigskip

\begin{abstract}
In  this paper we study two operations, Pasting and Reversing,
defined from a natural way to be applied over some rings such as
the ring of polynomials and the ring of linear differential
operators, which is a differential ring. We obtain some properties
of these operations over these rings, in particular over the set
of natural numbers, in where we rewrite some properties incoming
from recreational mathematics.\\

\footnotesize{\noindent{\bf Keywords and Phrases.} Differential rings, operators, pasting, polynomials, reversing}
\end{abstract}

\section{Introduction}
In this section we set the historical online and the theoretical
background  necessaries to understand the rest of the paper.\\

\subsection{Historical background}
Pasting and Reversing are common process for people any time. One
idea to become as mathematical operations these process was worked
for the first author since 1992. In particular, in 1993 were given
some lectures about the case of the natural numbers. These
lectures were published ten years after in \cite{Ac2}.\\

After, in
2008, these Pasting and Reversing operation were applied to obtain
families of {\it Simple permutations}, see  \cite{Ac1,AM}. \\

In the present work we introduce Pasting and Reversing operations
for the case of the ring of polynomials. In natural way is introduced Pasting and Reversing of natural numbers where some beautiful properties that appears in
\cite{Ta} are expressed
in terms of these operations. \\

Finally, we start the study of some properties of differential
rings. In particular from {\it differential Galois theory} (see
\cite{VS}), we start the analysis of these operations over linear
differential operators.

\subsection{Preliminaries}
We recall that a ring is a set $R$ equipped with two binary
operations $+ : R \times R \rightarrow R$ and $\cdot : R \times R
\rightarrow R$ satisfying the following requirements:
\begin{enumerate}
\item $\forall a, b \in R$, $a + b \in R$.
\item $\forall a, b, c \in R$, $a + (b + c)= (a + b) + c$.
\item $\exists 0 \in R$, such that $\forall a \in R$, $0 + a = a + 0 = a$.
\item $\forall a \in R$, $\exists b\in R$ such that $a + b = b + a =
0$
\item $\forall a, b \in R$, $a + b = b + a$.
\item $\forall a, b \in R$, $a \cdot b\in R$.
\item $\forall a, b, c\in R$, $(a \cdot b) \cdot c = a \cdot (b \cdot c)$ .
\item $\exists 1\in R$, such that $\forall a\in R$, $1 \cdot a = a \cdot 1 = a$.
\item $\forall a, b, c\in R$, $a \cdot (b + c) = (a \cdot b) + (a \cdot c)$.
\item $\forall a, b, c\in R$, $(a + b) \cdot c = (a \cdot c) +
(b \cdot c)$.
\end{enumerate} This definition can be found in any book of algebra see for example
\cite{He}.\\

We say that $R$ is a differential ring if  there exists a
derivation $\partial$ such that $$ \forall a, b\in R,
\partial(a+b)=\partial a + \partial b, \quad \partial (a \cdot
b)= \partial a \cdot b + a \cdot \partial b,$$ for more details see \cite{VS}. \\

In particular, we are interested in the ring of polynomials
$\mathbb{C} [x]$ and in linear differential operators
$\mathcal{L}$ defined as $$\mathcal{L}:= \sum^{n}_{k=0} a_{k}
\partial ^{k}, \quad \textrm{where}\quad a_{k} \in K,$$ being $K$ a differential field. The set of operators
$\mathcal{L}$  is a differential ring.
\medskip
\section{Pasting and Reversing over Polynomials}

In this section we introduce the Pasting and Reversing operations
over the ring of polynomials and over the set of natural numbers
in where some aspects of recreational mathematics are shown.\\

\subsection{Polynomial case}
We only consider polynomials $P$ such that $x \nmid P(x)$. For
suitability, we write $P$ as follows:$$ P(x)= \sum^{n}_{k=0}
a_{n-k}x^{n-k}.$$ We can see clearly that $1+{\rm deg} (P)= \c{C}
(P)$ is the number of coefficients of the polynomial $P$.

\begin{definition} [Reversing of Polynomials] \label{def1}
Consider $P \in \mathbb{C}[x]$ written as
\begin{equation}\label{eq1}
P(x)=\sum_{k=0}^{n} a_{n-k}x^{n-k},
\end{equation}
the Reversing of $P$, denoted by $\widetilde{P}$ is given by
\begin{equation}\label{eq2}
\widetilde{P}(x)=\sum_{k=0}^{n} b_{n-k}x^{n-k}, \quad
b_{n-k}=a_{k}, \quad k=0,1, \ldots, n.
\end{equation}
\end{definition}
\medskip
Although to work with ${\rm deg}$ is equivalent to work with
\c{C}, we prefer the last one for our convenience.\\

Definition \ref{def1} lead us to the following result.
\begin{proposition} \label{prop1}
Consider the polynomials $P$ and $\widetilde{P}$ as in equations
\eqref{eq1}, \eqref{eq2}  respectively. The following statement
holds.
\begin{enumerate}
\item $\widetilde{P}(x)=x^{n}P(1/x)$, where $n+1=\c{C} (P)$.
\item $\widetilde{P}(1/\alpha)=0$ if and only if $P(\alpha)=0$.
\item $\widetilde{P}(x)=(-1)^n(\alpha_1x-\beta_1)(\alpha_2x-\beta_2)\ldots
(\alpha_nx-\beta_n)$ if and only if $P(x)=(\beta_1x-\alpha_1)\cdots
(\beta_nx-\alpha_n)$.

\item $\widetilde{\widetilde{P}}=P$.
\item $\c{C}(P)=\c{C}(\widetilde{P})$.
\item $\widetilde{P+Q}=\widetilde{P}+\widetilde{Q}$, for
$\c{C}(P)=\c{C}(Q)$.
\item $\widetilde{P\cdot Q}=\widetilde{P}\cdot\widetilde{Q}$.
\end{enumerate}
\end{proposition}
\medskip \textbf{Proof.}  We start the proof according to each item:
\begin{enumerate}
\item By equation \eqref{eq1}, we can see that $$x^{n}P(1/x)=x^{n}\sum_{k=0}^{n}
a_{n-k}(1/x)^{n-k},$$so that $$ x^{n}P(1/x)=\sum_{k=0}^{n}
a_{n-k}x^{k}=\sum_{k=0}^{n} a_{k}x^{n-k},$$ In this way, by
equation \eqref{eq2}, we have $x^{n}P(1/x)=\widetilde{P}(x).$
\item Due to $ x\nmid P(x)$, $\alpha\neq0$. Now, by item $1$, taking
$x=(1/\alpha)$, we have
$$\widetilde{P}(1/\alpha)=(1/\alpha)^{n}P(\alpha).$$ By hypothesis $P(\alpha)=0$, for instance $\widetilde{P}(1/\alpha)=0$. In similar way for the converse.
\item From item $1$, we have $$\widetilde{P}(x)=x^{n}(\beta_{1}(1/x)- \alpha_{1})\cdots(\beta_{n}(1/x)- \alpha_{n})$$ in this way, $$\widetilde{P}(x)=(\beta_{1}- \alpha_{1}x)\cdots(\beta_{n}-
\alpha_{n}x),$$ for instance
$$\widetilde{P}(x)=(-1)^{n}(\alpha_{1}x-
\beta_{1})\cdots(\alpha_{n}x- \beta_{n}).$$ In similar way for the converse.
\item Assume $P(x)$ and $\widetilde{P}(x)$ as in  item $3$. Thus, we have
$$\widetilde{\widetilde{P}}(x)=(-x)^{n}(\alpha_{1}(1/x)-\beta_{1})\cdots(\alpha_{n}(1/x)-\beta_{n}),$$
so that
$$\widetilde{\widetilde{P}}(x)=(\beta_{1}x-\alpha_{1})\cdots(\beta_{n}x-\alpha_{n}),$$
for instance $\widetilde{\widetilde{P}}(x)=P(x)$.
\item From item $3$ we observe that ${\rm deg}(P)={\rm deg}(\widetilde{P}),$ thus $\c{C}(P)=\c{C}({\widetilde{P}}).$
\item Assume that $$P(x)=\sum^{n}_{k=0}a_{n-k}x^{n-k},\,
Q(x)=\sum^{n}_{k=0}b_{n-k}x^{n-k}.$$ Setting $R=P+Q$ we have that
$$R(x)=\sum^{n}_{k=0}c_{n-k}x^{n-k}, \quad c_j = a_j +
b_j,\quad j=0,\cdots,n.$$ Now, by equation \eqref{eq2} it follows
that
$$\widetilde{R}(x)=\sum^{n}_{k=0}c_{k}x^{n-k}=\sum^{n}_{k=0}a_{k}x^{n-k}+\sum^{n}_{k=0}b_{k}x^{n-k},$$
which means, again by equation \eqref{eq2}, that
$\widetilde{R}=\widetilde{P}+\widetilde{Q}$. Thus, we conclude
that $\widetilde{P+Q}=\widetilde{P}+\widetilde{Q}.$
\item Assume that $$P(x)=(\beta_{1}x-\alpha_{1})\cdots(\beta_{n}x-\alpha_{n}),\,
Q(x)=(\gamma_{1}x-\mu_{1})\cdots(\gamma_{m}x-\mu_{m}).$$ Setting
$R=P\cdot Q$ we have that
$$R(x)=(\beta_{1}x-\alpha_{1})\cdots(\beta_{n}x-\alpha_{n})(\gamma_{1}x-\mu_{1})\cdots(\gamma_{m}x-\mu_{m}).$$ By item
$3$ of Proposition \eqref{prop1} we have that
$$\widetilde{R}(x)=(-1)^{n+m}(\alpha_{1}x-\beta_{1})\cdots(\beta_{n}x-\alpha_{n})(\gamma_{1}x-\mu_{1})\cdots(\gamma_{m}x-\mu_{m}),$$
being $\widetilde{R}(x)$ equivalent to
$$[(-1)^{n}(\alpha_{1}x-\beta_{1})\cdots(\beta_{n}x-\alpha_{n})][(-1)^m(\gamma_{1}x-\mu_{1})\cdots(\gamma_{m}x-\mu_{m})].$$
Thus, we obtain
$\widetilde{R}(x)=\widetilde{P}(x)\cdot\widetilde{Q}(x)$, which
implies that $\widetilde{P\cdot Q
}=\widetilde{P}\cdot\widetilde{Q}$.
\end{enumerate}

\begin{remark}
From the start we assumed that $x\nmid P(x)$. In case that $x\mid
P(x)$, items $2$ and $5$ are false. We recall that items $4$ and $5$ also can be proven using only Definition \ref{def1}, i.e., equations \eqref{eq1} and \eqref{eq2}.
\end{remark}

There are some specific known cases in which we can use the
reversing operation over special families of polynomials such as
\emph{Bessel polynomials}, see \cite{Br,Gr}.

Definition \ref{def1} lead us to the following definition.
\begin{definition} \label{def2} Polynomials $P$ and $Q$
are called palindromic and antipalindromic polynomials
respectively whether they satisfy
$$\widetilde{P}=P,\quad \widetilde{Q}=-Q.$$

\end{definition}

Proposition \ref{prop1} and Definition \ref{def2}, lead us to the
following results.
\begin{proposition} \label{prop2}
Let $P$ be a palindromic or antipalindromic polynomial
with roots $\alpha_{1},\cdots,\alpha_{n}$, being $n+1=\c{C}(P)$, then

 $$\alpha_{k+j}=1/\alpha_ j,\quad \c{C}(P)\in\{2k+1,2k+2\},\quad
 j=1,\cdots,k.$$ Furthermore, if $\c{C}(P)=2k+2$ and $P$ is palindromic (respectively antipalindromic),
 then $\alpha_{2k+1}=-1$ (respectively $\alpha_{2k+1}=1$).

\end{proposition}
\medskip \textbf{Proof.}
By Definition \ref{def2} and Proposition \ref{prop1},
$P(\alpha_{j})=0$ implies that
$$\widetilde{P}(1/\alpha_{j})=\pm{P}(1/\alpha_{j})=0,\quad
j=1,\cdots,\c{C}(P)-1.$$ Thus, for $\c{C}(P)=2k+1$ we can arrange
$\alpha_{k+j}=1/\alpha_{j}$, $j=1,\cdots,k$. In the same  way for
$\c{C}(P)=2k+2$, we have $\alpha_{k+j}=1/\alpha_{j}$,
$j=1,\cdots,k$ and $\alpha_{2k+1}=1/\alpha_{2k+1}$, so that
$\alpha_{2k+1}=\pm1$. If $P=\widetilde{P}$, then the signs of the
its coefficients must be preserved, so that $\alpha_{2k+1}$ must
be $-1$. Finally, if $\widetilde{P}=-P$ then the signs of the
coefficients must be interchanged, so that $\alpha_{2k+1}$ must be
$1$.

\begin{proposition}\label{prop3}
The following statements holds.
\begin{enumerate}
\item The addition of two palindromic polynomials, with the same degree, is also a
palindromic polynomial.
\item The product of two palindromic polynomials is also a
palindromic polynomial.
\item The addition of two antipalindromic polynomials, with the same degree, is also an
antipalindromic polynomial.
\item The product of two antipalindromic polynomials is a
palindromic polynomial.
\item The product of a palindromic polynomial with an antipalindromic polynomial is
an antipalindromic polynomial.
\end{enumerate}
\end{proposition}

\medskip \textbf{Proof.}  We prove the proposition according each item.
\begin{enumerate}
\item Let $P$ and $Q$ be palindromic polynomials. By item $6$ of
Proposition \ref{prop1} we have that
$\widetilde{P+Q}=\widetilde{P}+\widetilde{Q}=P+Q$. In consequence,
$P+Q$ is a palindromic polynomial.
\item Let $P$ and $Q$ be palindromic polynomials. By item $7$ of
Proposition \ref{prop1} we have that $\widetilde{P\cdot
Q}=\widetilde{P}\cdot\widetilde{Q}=P\cdot Q$. In consequence,
$P\cdot Q$ is a palindromic polynomial.
\item Let $P$ and $Q$ be antipalindromic polynomials. By item $6$ of
Proposition \ref{prop1} we have that
$\widetilde{P+Q}=\widetilde{P}+\widetilde{Q}=-P-Q=-(P+Q)$. In
consequence, $P+Q$ is an antipalindromic polynomial.
\item Let $P$ and $Q$ be antipalindromic polynomials. By item $7$ of
Proposition \ref{prop1} we have that $\widetilde{P\cdot
Q}=\widetilde{P}\cdot\widetilde{Q}=(-P)\cdot (-Q)=PQ$. In
consequence, $P\cdot Q$ is a palindromic polynomial.
\item Let $P$ be a palindromic polynomial and let be $Q$ an antipalindromic polynomial. By item $7$ of
Proposition \ref{prop1} we have that $\widetilde{P\cdot
Q}=\widetilde{P}\cdot\widetilde{Q}=P\cdot (-Q)=-P\cdot Q$. In
consequence, $P\cdot Q$ is an antipalindromic polynomial.

\end{enumerate}

The following definition corresponds to a natural example of \emph{orthogonal polynomials}, see \cite{ch,Ma,Ri}.
\begin{definition}[Chebyshev polynomials of first kind]\label{chedef} The Chebyshev polynomials of the first kind, denoted by $T_n$, are defined by the trigonometric identity
\[T_n(w)=\cos(n \arccos w)=\cosh(n\,\mathrm{arccosh}\,w),\quad n\in\Z_+,\]

which is equivalent to the identity

\[T_n(\cos(\alpha))=\cos(n\alpha)=\cosh(n\alpha).\]

\end{definition}
\begin{lemma}\label{the1}
If $w$ is given by $\frac{1}{2}(z+\frac{1}{z})$ then $\frac{1}{2}(z^{n}+\frac{1}{z^{n}})=T_{n}(w)$
\end{lemma}
\medskip \textbf{Proof.} By Definition \ref{chedef} we write
$$T_{1}(w)=\cos(\alpha)=w$$
$$\vdots$$
$$T_{n}(w)=\cos(n\alpha),$$ which lead us to $$T_{n}(w)=\frac{e^{in\alpha}+e^{-in\alpha}}{2}=\frac{\left(e^{i\alpha}\right)^{n}+\left(e^{i\alpha}\right)^{-n}}{2}=\frac{1}{2}(z^{n}+z^{-n}).$$
In particular, $w=\frac{1}{2}(z+\frac{1}{z})$.\medskip

The following result has been suggested by V. Sokolov.

\begin{proposition}\label{sokoth}
Let $P_{2n}$ be a palindromic polynomial with coefficients $a_i$, $0\leq i\leq 2n$. Then
\begin{equation}\label{sokoeq}\frac{P_{2n}(z)}{2z^{n}}=\sum_{k=0}^{n}a_{n-k}T_{k}(w),\quad w=\frac{1}{2}\left(z+\frac{1}{z}\right).\end{equation}
\end{proposition}
\medskip \textbf{Proof.} By hypothesis $P_{2n}(z)=a_{2n}z^{2n}+a_{2n-1}z^{2n-1}+\ldots +a_{1}z+a_{0}$. Now, due to $P_{2n}=\widetilde{P}_{2n}$, we have that $a_{i}=a_{2n-i}$, $i=0,\ldots,2n.$ Thus, dividing $P_{2n}(z)$ by $2z^n$ we obtain $$\frac{P_{2n}(z)}{2z^{n}}=\frac{a_{2n}z^{n}+a_{2n-1}z^{n-1}+\ldots +a_{2n-1}z^{1-n}+a_{2n}z^{-n}}{2}.$$ Reorganizing the common coefficients we have  $$\frac{P_{2n}(z)}{2z^{n}}=\frac{a_{2n}}{2}\left(z^{n}+\frac{1}{z^{n}}\right)+\frac{a_{2n-1}}{2}\left(z^{n-1}+\frac{1}{z^{n-1}}\right)+\ldots+\frac{1}{2}a_{n}.$$ By the change of variable $w=\frac{1}{2}(z+\frac{1}{z})$ and Lemma \ref{the1}, the righthand side is
$a_{2n}T_n(w)+a_{2n-1}T_{n-1}(w)+\ldots+a_nT_0(w)$ when $a_i=a_{2n-i}$, $i=0,..,2n$. Thus we can conclude the expression given in equation \eqref{sokoeq}.\medskip

 The concept of palindromic and antipalindromic polynomials is
 very ancient, there are a lot of references about these
 polynomials using the concept of \emph{Reciprocal Polynomials}, see for example
 \cite{Bo,Ev,Ghm,Hi,La,Sa} and references therein. We recall, using the previous references, that $P$ is the reciprocal of $Q$ if $Q=\widetilde{\overline
{P}}=P^{\ast}, \quad z=a+bi, \quad \overline{z}=a-bi$,
$$P(z)=a_{n}z^{n}+a_{n-1}z^{n-1}+\ldots+a_{1}z+a_{0}$$
$$Q(z)=\widetilde{\overline{P}}(z)=a_{0}\overline{z^{n}}+a_{1}\overline{z^{n-1}}+\ldots+a_{n-1}\overline{z}+a_{n}.$$
On another hand, $P$ is self reciprocal if
$P=P^{\ast}=\widetilde{\overline P}$. It is easy to see that for
$b=0$, this means that $z\in\mathbb{R}$, then $P=\widetilde{P}$,
that is, $P$ is a palindromic polynomial. In the same way, for
antipalindromic polynomials.\\

Now we introduce the definition of Pasting operation over
polynomials.

\begin{definition}\label{def3} Pasting of the polynomials $P$ and $Q$, denoted by $P\diamond Q$, is given by:
$x^{\c{C}(Q)}P+Q$.
\end{definition}

The following properties are consequences of Definition
\ref{def3}.

\begin{proposition}\label{prop4} Let $P, Q, R$ be polynomials. The following statements holds:
\begin{enumerate}
\item $\tilde{P}\diamond \tilde{Q}=\widetilde{{Q\diamond P}}$
\item $(P\diamond Q)\diamond R= P\diamond (Q\diamond R)$
\end{enumerate}
\end{proposition}
\medskip \textbf{Proof.} We consider separately each item.
\begin{enumerate}
\item Let $P, Q$ be polynomials, where $$P(x)=\sum_{l=0}^{s} a_{s-l}x^{s-l},\quad Q(x)= \sum_{j=0}^{k}
b_{k-j}x^{k-j}.$$ Then, by Definition \ref{def3} and assuming
$R=Q\diamond P$ we have
$$R(x)=x^{s+1}\sum_{j=0}^{k}b_{k-j}x^{k-j}+
\sum_{l=0}^{s} a_{s-l}x^{s-l}=\sum_{i=0}^{k+s+1}
c_{k+s+1-i}x^{k+s+1-i},$$ where the coefficients $c_m$ are given
by
$$c_m= \left\{\begin{array}{l}
  a_m,\quad 0\leq m\leq s,  \\
  b_m, \quad s+1 \leq m \leq
k+s+1.
\end{array}\right.$$

By Definition \ref{def1} we obtain
$$\widetilde{P}(x)=\sum_{l=0}^{s} a_{l}x^{s-l},\,\, \widetilde{Q}(x)= \sum_{j=0}^{k}
b_{j}x^{k-j},\,\, \widetilde{R}(x)=\sum_{i=0}^{k+s+1}
c_{i}x^{k+s+1-i}.$$ In consequence, by Definition \ref{def3} we
have that $$\widetilde{R}(x)=\sum_{i=0}^{k+s+1}
c_{i}x^{k+s+1-i}=x^{k+1}\sum_{l=0}^{s}a_{l}x^{s-l}+\sum_{j=0}^{k}
b_{j}x^{k-j}.$$ So that we obtain
$\widetilde{R}=\widetilde{P}\diamond \widetilde{Q}.$

\item Let $P, Q, R$ polynomials, such that $$P(x)=\sum_{i=0}^{k}a_{k-i}x^{k-i},\,\, Q(x)=\sum_{i=0}^{j}
b_{j-i}x^{j-i},\,\, R(x)=\sum_{i=0}^{l} d_{l-i}x^{l-i},$$ where,
$\c{C}(P)=k+1, \c{C}(Q)=j+1$ and $\c{C}(R)=l+1$. We write $(P \diamond
Q)\diamond R$ by means of Definition \ref{def1} as follows
$$x^{l+1}\left(x^{j+1}\sum_{i=0}^{k}a_{k-i}x^{k-i}+\sum_{i=0}^{j} b_{j-i}x^{j-i}\right)+ \sum_{i=0}^{l}
d_{l-i}x^{l-i},$$ it is rewritten as
$$x^{j+l+2}\sum_{i=0}^{k}a_{k-i}x^{k-i}+ \left(x^{l+1}\sum_{i=0}^{j} b_{j-i}x^{j-i}+
\sum_{i=0}^{l} d_{l-i}x^{l-i}\right),$$ which can be written as
$$x^{j+l+2}\sum_{i=0}^{k}
a_{k-i}x^{k-i}+ \left(\sum_{i=0}^{j} b_{j-i}x^{j-i}\diamond
\sum_{i=0}^{l} d_{l-i}x^{l-i}\right).$$ In consequence $$(P \diamond
Q)\diamond R=P \diamond(Q\diamond R).$$
\end{enumerate}

As consequence of Proposition \ref{prop2}, which is adapted for
pasting of polynomials, we present the following result.
\begin{proposition}\label{prop5} Let $P$ be a polynomial. The linear polynomial $x+1$ divides to the polynomial $P\diamond\widetilde{P}$.
\end{proposition}
\medskip \textbf{Proof.}  Owing to $\c{C}(P\diamond\widetilde{P})$ is even and $P\diamond\widetilde{P}$
is palindromic, then, by Proposition \ref{prop2}, $-1$ is root of
$P\diamond\widetilde{P}$, so that $x+1\mid P\diamond\widetilde{P}$.

\subsection{Natural numbers case}
\medskip
The properties presented before in the polynomial case are very
useful for natural numbers choosing $x=10$. For natural numbers,
$\c{C}$ is called \emph{digital cipher}, see \cite{Ac2}. \\

We recall that, by previous results, the reversing of
$n\in\mathbb{N}$ is given by
$$\widetilde{n}=\sum_{j=0}^{r} a_{j}10^{r-j},\quad \textrm{where}\quad n=\sum_{j=0}^{r}
a_{r-j}10^{r-j}.$$In a natural way, we introduce the concept of
\emph{palindrome numbers}: $n$ in palindrome if and only if
$n=\widetilde{n}.$ In the same way, the pasting of $n,\, m \in
\mathbb{N}$ is given by $10^{\c{C}(m)}n +m$.
\\
For natural number case Propositions \ref{prop1}, \ref{prop4} and
\ref{prop5}, can be summarized in the following result.

\begin{proposition} Let $n,\, m,\, p \in\mathbb{N}$, the
following statements holds:
\begin{enumerate}
\item $\widetilde{\tilde{n}}=n$
\item $\tilde{n}\diamond \tilde{m}=\widetilde{{m\diamond n}}$
\item $(m\diamond n)\diamond p= m\diamond (n\diamond p)$
\item If $n$ is palindrome and $\c{C}({n})$ is even, then $11$ is a
divisor of $n$.
\item $11|n \diamond \tilde{n}$.
\end{enumerate}
\end{proposition}
In general the properties presented in polynomial case for the
operations of the ring ($+,\cdot$) are not true for natural
numbers, although doing some restrictions we can obtain similar
results as presented before.\medskip

As application of pasting and reversing operations over natural
numbers, we can rewrite some mathematical games such as the
presented in \cite{Ta}. For suitability to our purposes, we
introduce the following notation
\begin{equation}\label{eq3}
\diamondsuit_{k=0}^n a_k:=a_0\diamond a_1\diamond \ldots\diamond
a_n.
\end{equation}
The following mathematical games can be found in \cite{Ta}, but here we use our
approach.
\begin{enumerate}
\item $(\diamondsuit_{k=0}^n9-k)\cdot 9+9-(n+2)=\diamondsuit_{k=0}^{n+1}8$,
where $n\leq 9$. This can de expanded as:
\medskip

$\begin{array}{c}9\times 9+7=88\\ 98\times 9+6=888\\987\times
9+5=8888\\9876\times 9+4=88888\\98765\times
9+3=888888\\987654\times 9+2=8888888\\9876543\times
9+1=88888888\\98765432\times 9+0=888888888\\987654321\times
9-1=8888888888
\end{array}$
\medskip

\item We know that $1^2=1$, now, for $0< n<9$,we have
$$(\diamondsuit_{k=0}^{n}1)^2=\diamondsuit_{k=0}^{n}(k+1)\diamond\widetilde{\diamondsuit_{k=0}^{n-1}(k+1)}.$$ This can be expanded
as:
\medskip

$\begin{array}{c}1\times 1=1\\ 11\times 11=121\\ 111\times
111=12321\\ 1111\times 1111=1234321\\ 11111\times 11111=123454321\\
111111\times 111111=12345654321\\ 1111111\times
1111111=1234567654321\\ 11111111\times 11111111=123456787654321\\
111111111\times 111111111=12345678987654321
\end{array}$
\end{enumerate}

\section{Pasting and Reversing over Differential Operators}

We consider linear differential operators $$\mathcal L=a_n\partial^n+a_{n-1}\partial^{n-1}+\ldots+a_1\partial+a_0,\quad a_0\neq 0,\quad a_i\in K,$$ where $i=0,1,\ldots,n$ and $K$ is a \textit{differential field}, see \cite{VS}. Solutions of linear differential equations are related with the factorization of linear differential operators, see \cite{Ho}.\medskip

From now on we understand as differential operators the linear differential operators. For
suitability, we write $\mathcal L$ as follows:$$ \mathcal L= \sum^{n}_{k=0}
a_{n-k}\partial^{n-k}.$$ As in previous cases, we denote by $\c{C}
(\mathcal L)$ the number of coefficients of the differential operator $\mathcal L$. For instance, if the order of $\mathcal L$ is $n$, then $\c{C}
(\mathcal L)=n+1$.

\begin{definition} [Reversing of Differential Operators] \label{defdo1}
Let\\ $\mathcal L$ be a differential operator written as
\begin{equation}\label{eqdo1}
\mathcal L=\sum_{k=0}^{n} a_{n-k}\partial^{n-k}, \quad a_0\neq 0,\quad a_i\in K.
\end{equation}
The Reversing of $\mathcal L$, denoted by $\widetilde{\mathcal L}$ is given by
\begin{equation}\label{eqdo2}
\widetilde{\mathcal L}=\sum_{k=0}^{n} b_{n-k}\partial^{n-k}, \quad
b_{n-k}=a_{k}, \quad k=0,1, \ldots, n.
\end{equation}
\end{definition}
\medskip

Definition \ref{defdo1} lead us to the following result.
\begin{proposition} \label{propdo1}
Let $\mathcal L$ and $\widetilde{\mathcal L}$ be differential operators as in equations
\eqref{eqdo1}, \eqref{eqdo2}  respectively. The following statement
holds.
\begin{enumerate}
\item $\widetilde{\widetilde{\mathcal L}}=\mathcal L$.
\item $\c{C}(\mathcal L)=\c{C}(\widetilde{\mathcal L})$.
\item $\widetilde{\mathcal L+ \mathcal R}=\widetilde{\mathcal L}+\widetilde{\mathcal R}$, for
$\c{C}(\mathcal L)=\c{C}(\mathcal R)$.
\end{enumerate}
\end{proposition}
\medskip \textbf{Proof.}  Items $1$ and $2$ are consequences of the Definition \ref{defdo1}. Item $3$ is proven in similar way to polynomial case.\medskip

\begin{remark}
In general, there is not relationship between  $\ker \mathcal L$ and $\ker\mathcal{\widetilde L}$. Using \emph{Differential Galois Theory}, see \cite{VS}, it can be shown that $e^{-\frac{x^2}2}\in\ker (\partial^2+1-x^2)$, while there are not \emph{Liouvillian functions} in $\ker ((1-x^2)\partial^2+1)$.
\end{remark}

As particular case, we have the following result.
\begin{proposition}
Assume $\c{C} (\mathcal L)=2$, $y\in \ker \mathcal L$ and $u\in \ker \widetilde{\mathcal L}$. Then $(\partial\ln y)(\partial\ln u)=1$.
\end{proposition}

\noindent \textbf{Proof.}  Solving the linear differential equations $\mathcal L y=0$ and $\widetilde{\mathcal L} u=0$, we obtain $$e^{-\int\frac{a_0}{a_1}}\in\ker \mathcal L,\quad e^{-\int\frac{a_1}{a_0}}\in\ker \widetilde{\mathcal L}.$$ Thus, $(\partial\ln y)(\partial\ln u)=1$.\medskip

Definition \ref{defdo1}, as in polynomial case, lead us to the following definition.
\begin{definition} \label{defdo2} Differential operators $\mathcal L$ and $\mathcal R$
are called palindromic and antipalindromic differential operators
respectively whether they satisfy
$$\widetilde{\mathcal L}=\mathcal L,\quad \widetilde{\mathcal R}=-\mathcal R.$$

\end{definition}

Proposition \ref{propdo1} and Definition \ref{defdo2} lead us to the
following results.
\begin{proposition} \label{propdo2}
Let $\mathcal L$ and $\mathcal R$ be palindromic and antipalindromic differential operators respectively, being $\c{C}(\mathcal L)=\c{C}(\mathcal R)=2k$. Then there exist differential operators $\mathcal S$ and $\mathcal T$ such that $$\mathcal L=\mathcal S(\partial+1),\quad \mathcal R=\mathcal T(\partial-1).$$

\end{proposition}

\noindent\textbf{Proof.} We can see that $e^{-x}\in \ker(\partial +1)$ and $e^{x}\in \ker(\partial -1)$. Now, due to $\widetilde{\mathcal L}=\mathcal L$, then $a_{2k-1-i}=a_i$ and $\partial^{2k-1-i}e^{-x}=-\partial^{i}e^{-x}$. In this way, $e^{-x}\in \ker\mathcal L$, which means that there exists $\mathcal S$ such that $\mathcal L=\mathcal S(\partial + 1)$. On another hand, owing to $\widetilde{\mathcal R}=-\mathcal R$ then $a_{2k-1-i}=-a_i$ and $\partial^{2k-1-i}e^{x}=\partial^{i}e^{x}$. In this way, $e^{x}\in \ker\mathcal R$, which means that there exists $\mathcal T$ such that $\mathcal R=\mathcal T(\partial - 1)$.\medskip

We recall that differential operators $\mathcal T$ and $\mathcal S$ are \emph{left divisors} of $\mathcal L$ and $\mathcal R$ respectively. In the same way, $\partial+1$ and $\partial-1$ are \emph{right divisors} of $\mathcal L$ and $\mathcal R$ respectively. For further details see \cite{Ho}.

\begin{proposition}\label{propdo3}
The following statements holds.
\begin{enumerate}
\item The addition of two palindromic differential operators, with the same order, is also a
palindromic differential operator.

\item The addition of two antipalindromic differential operators, with the same order, is also an
antipalindromic differential operator.
\end{enumerate}
\end{proposition}

\noindent \textbf{Proof.} We proceed exactly as in Proposition \ref{prop3} for the polynomial case, using Definition \ref{defdo2} and item $3$ of Proposition \ref{propdo1}.\medskip

Now we introduce the definition of Pasting operation over
differential operators.

\begin{definition}\label{defdo3} Pasting of the differential operators $\mathcal{L}$ and $\mathcal R$, denoted by $\mathcal L\diamond \mathcal R$, is given by:
$\mathcal L\partial^{\c{C}(\mathcal R)}+\mathcal L$.
\end{definition}

The following properties, adapted from Proposition \ref{prop4}, are consequences of Definition
\ref{defdo3}.

\begin{proposition}\label{propdo4} Let $\mathcal L$, $\mathcal R$ and $\mathcal S$ be differential operators. The following statements holds:
\begin{enumerate}
\item $\tilde{\mathcal L}\diamond \tilde{\mathcal R}=\widetilde{{\mathcal R\diamond \mathcal L}}$
\item $(\mathcal L\diamond \mathcal R)\diamond \mathcal S= \mathcal L\diamond (\mathcal R\diamond \mathcal S)$
\end{enumerate}
\end{proposition}
\medskip \textbf{Proof.} We consider separately each item.
\begin{enumerate}
\item Let $\mathcal L, \mathcal R$ be differential operators, where $$\mathcal L=\sum_{l=0}^{s} a_{s-l}\partial^{s-l},\quad \mathcal R= \sum_{j=0}^{k}
b_{k-j}\partial^{k-j}.$$ Then, by Definition \ref{defdo3} and assuming
$\mathcal S=\mathcal R\diamond \mathcal L$ we have
$$\mathcal S=\left(\sum_{j=0}^{k}b_{k-j}\partial^{k-j}\right)\partial^{s+1}+
\sum_{l=0}^{s} a_{s-l}\partial^{s-l}=\sum_{i=0}^{k+s+1}
c_{k+s+1-i}\partial^{k+s+1-i},$$ where the coefficients $c_m$ are given
by
$$c_m= \left\{\begin{array}{l}
  a_m,\quad 0\leq m\leq s,  \\
  b_m, \quad s+1 \leq m \leq
k+s+1.
\end{array}\right.$$

By Definition \ref{defdo1} we obtain
$$\widetilde{\mathcal L}=\sum_{l=0}^{s} a_{l}\partial^{s-l},\,\, \widetilde{\mathcal R}= \sum_{j=0}^{k}
b_{j}\partial^{k-j},\,\, \widetilde{\mathcal S}=\sum_{i=0}^{k+s+1}
c_{i}\partial^{k+s+1-i}.$$ In consequence, by Definition \ref{defdo3} we
have that $$\widetilde{\mathcal S}=\sum_{i=0}^{k+s+1}
c_{i}\partial^{k+s+1-i}=\left(\sum_{l=0}^{s}a_{l}\partial^{s-l}\right)\partial^{k+1}+\sum_{j=0}^{k}
b_{j}\partial^{k-j}.$$ So that we obtain
$\widetilde{\mathcal S}=\widetilde{\mathcal L}\diamond \widetilde{\mathcal R}.$

\item Let $\mathcal L, \mathcal R, \mathcal S$ differential operators, such that $$\mathcal L=\sum_{i=0}^{k}a_{k-i}\partial^{k-i},\,\, \mathcal R=\sum_{i=0}^{j}
b_{j-i}\partial^{j-i},\,\, \mathcal S=\sum_{i=0}^{l} d_{l-i}\partial^{l-i},$$ where,
$\c{C}(\mathcal L)=k+1, \c{C}(\mathcal R)=j+1$ and $\c{C}(\mathcal S)=l+1$. We write $(\mathcal L \diamond
\mathcal R)\diamond \mathcal S$ by means of Definition \ref{defdo1} as follows
$$\left(\left(\sum_{i=0}^{k}a_{k-i}\partial^{k-i}\right)\partial^{j+1}+\sum_{i=0}^{j} b_{j-i}\partial^{j-i}\right)\partial^{l+1}+ \sum_{i=0}^{l}
d_{l-i}\partial^{l-i},$$ it is rewritten as
$$\left(\sum_{i=0}^{k}a_{k-i}\partial^{k-i}\right)\partial^{j+l+2}+ \left(\left(\sum_{i=0}^{j} b_{j-i}\partial^{j-i}\right)\partial^{l+1}+
\sum_{i=0}^{l} d_{l-i}\partial^{l-i}\right),$$ which can be written as
$$\left(\sum_{i=0}^{k}
a_{k-i}\partial^{k-i}\right)\partial^{j+l+2}+ \left(\sum_{i=0}^{j} b_{j-i}\partial^{j-i}\diamond
\sum_{i=0}^{l} d_{l-i}\partial^{l-i}\right).$$ In consequence $$(\mathcal L \diamond
\mathcal R)\diamond \mathcal S= \mathcal L\diamond(\mathcal R\diamond \mathcal S).$$
\end{enumerate}

As consequence of Proposition \ref{propdo2}, which is adapted for
pasting of polynomials, we present the following result.

\begin{proposition}\label{propdo5} Let $\mathcal L$ be a differential operator. The operator $\partial+1$ is a right divisor of the differential operator $\mathcal L\diamond\widetilde{\mathcal L}$.
\end{proposition}
\medskip \textbf{Proof.}  Owing to $\c{C}(\mathcal L\diamond\widetilde{\mathcal L})$ is even and $\mathcal L\diamond\widetilde{\mathcal L}$
is palindromic, then, by Proposition \ref{propdo2}, $e^{-x}\in\ker P\diamond\widetilde{P}$, so that $\mathcal T(\partial+1)=P\diamond\widetilde{P}$, for some differential operator $\mathcal T$.\medskip

The following results are particular cases of differential operators in where the differential field is considered as $K=\mathbb{C}$.

\begin{proposition} \label{propdoc1}
Consider $\mathcal L$ and $\widetilde{\mathcal L}$ as in equations
\eqref{eqdo1}, \eqref{eqdo2}  respectively, being $K=\mathbb{C}$. The following statement
holds.
\begin{enumerate}
\item $\widetilde{\mathcal L}=\partial^{n}\sum_{k=0}^na_{n-k}\partial^{k-n}$, where $\mathcal L=\sum_{k=0}^na_{n-k}\partial^{n-k}$.
\item $x^ke^{-\lambda_ix}\in\ker\widetilde{\mathcal L}$ if and only if $x^ke^{\lambda_ix}\in\ker\mathcal L$.
\item $\widetilde{\mathcal L}=(-1)^n(\alpha_1\partial-\beta_1)(\alpha_2\partial-\beta_2)\ldots
(\alpha_n\partial-\beta_n)$ if and only if $\mathcal L=(\beta_1\partial-\alpha_1)(\beta_2\partial-\alpha_2)\cdots
(\beta_n\partial-\alpha_n)$, $\alpha_i,\beta_i\in\mathbb{C}$.
\item $\widetilde{\mathcal L\cdot \mathcal R}=\widetilde{\mathcal L}\cdot\widetilde{\mathcal R}$.
\item If $\mathcal L$ is a palindromic (or antipalindromic) differential operator such that $\{x^{i_1}e^{\lambda_{1}x},\cdots,x^{i_r}e^{\lambda_{r}x}\}\subset\ker\mathcal L$, then

 $$e^{\lambda_{k+j}x}=e^{-\lambda_ jx},\quad r\in\{2k,2k+1\},\quad
 j=1,\cdots,k.$$
\item The product of two palindromic differential operators is also a
palindromic differential operator.
\item The product of two antipalindromic differential operators is a
palindromic differential operator.
\item The product of a palindromic differential operator with an antipalindromic differential operator is
an antipalindromic differential operator.

\end{enumerate}
\end{proposition}
\medskip

\noindent \textbf{Proof.} It follows since the characteristic polynomial satisfy the same properties (see Propositions \ref{prop1}, \ref{prop2}, \ref{prop3}) and due to $\partial a=a\partial$ for all $a\in\mathbb{C}$.

\section*{Final Remark}
This paper is one starting point to develop several research projects such as:

\begin{itemize}
\item Applications of Pasting and Reversing over vectorial spaces and matrices.\medskip

\item Applications of Pasting and Reversing over polynomials in several variables.\medskip

\item Applications of Pasting and Reversing over general differential operators.\medskip

\item Applications of Pasting and Reversing over general difference and $q$ -difference operators.\medskip

\item Applications of Pasting and Reversing over general simple permutations and combinatorial dynamics\medskip.

\item Applications of Pasting and Reversing in physics, particularly in \emph{supersymmetric quantum mechanics}.
\end{itemize}

There are papers in which this approach can be applied, see for example \cite{Bo,Ghm} for the polynomial case. We hope that the material presented here can be useful for the interested reader.

\section*{Acknowledgments}
The research of the first author has been supported by Universidad Sergio Arboleda and by the MCyT-FEDER Grant MTM2006-00478, Spanish Government. Part of this research was presented by the second and third author at XVII Colombian Congress of Mathematics, developed in Cali on August 2009. The authors thank to Jes\'us Hernando P\'erez, Vladimir Sokolov and David Bl\'azquez-Sanz by their useful comments and suggestions on this work.

\end{document}